\documentclass[a4paper]{article}

\usepackage{graphicx}
\usepackage{amssymb}
\usepackage{latexsym}
\usepackage{verbatim}
\usepackage{amsthm}

\theoremstyle{definition}

\newtheorem{definition}{Definition}[section]
\theoremstyle{plain}
\newtheorem{theorem}{Theorem}[section]
\newtheorem{corollary}{Corollary}[section]

\newtheorem{lemma}{Lemma}[section]
\newtheorem{proposition}{Proposition}[section]

\def\sqr#1#2{{\vcenter{\vbox{\hrule height.#2pt
\hbox{\vrule width.#2pt height#1pt \kern#1pt \vrule width.#2pt}
\hrule height.#2pt}}}}
\def\square{\mathchoice\sqr55\sqr55\sqr{2.1}3\sqr{1.5}3}

\title{\bf Duality index of oriented regular hypermaps}

\author{Daniel Pinto\\
\small CMUC, Department of Mathematics, University of Coimbra \\[-0.8ex]
\small 3001-454 Coimbra, Portugal \\
\small \texttt{dpinto@mat.uc.pt} }

\begin{document}

\maketitle

\begin{abstract}

By adapting the notion of chirality group, the duality group of
$\cal H$ can be defined as the the minimal subgroup $D({\cal H})
\trianglelefteq Mon({\cal H})$ such that ${\cal H}/D({\cal H})$ is
a self-dual hypermap (a hypermap isomorphic to its dual). Here, we
prove that for any positive integer $d$, we can find a hypermap of
that duality index (the order of $D({\cal H})$), even when some
restrictions apply, and also that, for any positive integer $k$,
we can find a non self-dual hypermap such that $|Mon({\cal
H})|/d=k$. This $k$ will be called the \emph{duality coindex} of
the hypermap.

\medskip



\end{abstract}

\section{Operations on hypermaps}

Topologically, a \emph{map} is a cellular embedding of a connected
graph into a closed connected surface. We can generalize this
notion if, instead of graphs, we use hypergraphs, allowing each
(hyper)edge to be adjacent to more than two (hyper)vertices. By
this process, we can construct a more general structure: the
\emph{hypermap}. Usually, hypermaps are represented by cellular
embeddings of bipartite maps (following the Walsh correspondence
between bipartite maps and hypermaps \cite{Walsh}) or by cellular
embeddings of connected trivalent graphs. In this last
representation (James representation \cite{James}) we label each
face 0, 1 or 2 so that each edge of the graph is incident to two
faces carrying different labels. Those faces correspond to
hypervertices, hyperedges or hyperfaces, depending on the label
they carry.

 The topological definitions of maps and
hypermaps have a combinatorial translation. A map can be
understood as a transitive permutation representation
$\Gamma\to{\rm Sym}\,F$ of the group
$$\Gamma=\langle r_0,r_1,r_2|r_0^2=r_1^2=r_2^2=(r_2r_0)^2=1
\rangle= V_4 * C_2$$ on a set $F$ representing its flags (the
cells of the barycentric subdivision of the map); and a hypermap
can be regarded as a transitive permutation representation
$\Delta\to{\rm Sym}\,\Omega$ of the group
$$\Delta=\langle r_0, r_1, r_2\mid r_0^2=r_1^2=r_2^2=1\rangle\cong
C_2*C_2*C_2,$$ on a set $\Omega$ representing its hyperflags.
Similarly, an oriented hypermap (without boundary) can be regarded
as a transitive permutation representation of the subgroup
$$\Delta^+=\langle \rho_0, \rho_1, \rho_2\mid \rho_0\rho_1\rho_2
=1\rangle=\langle \rho_0,\rho_2\mid -\rangle$$ of index $2$ in
$\Delta$ (a free group of rank $2$) consisting of the elements of
even word-length in the generators $r_i$, where $\rho_0=r_1r_2$,
$\rho_1=r_2r_0$ and $\rho_2=r_0r_1$. In the case of hypermaps, the
hypervertices, hyperedges and hyperfaces ($i$-dimensional
constituents for $i=0, 1, 2$) are the orbits of the dihedral
subgroups $\langle r_1, r_2\rangle$, $\langle r_2, r_0\rangle$ and
$\langle r_0, r_1\rangle$, and in the case of oriented hypermaps
they are the orbits of the cyclic subgroups
$\langle\rho_0\rangle$, $\langle\rho_1\rangle$ and
$\langle\rho_2\rangle$, with incidence given by nonempty
intersection in each case. The local orientation around each
hypervertex, hyperedge or hyperface is determined by the cyclic
order of the corresponding cycle of $\rho_0, \rho_1$ or $\rho_2$.

Operations on topological maps were first studied by S. Wilson
\cite{Operators} but were later extended to hypermaps, following a
more algebraic approach \cite{Thornton, Lynne}. If $\mathcal H$ is
a hypermap corresponding to a permutation representation
$\theta:\Delta\to{\rm Sym}\,\Omega$, and if $\alpha$ is an
automorphism of $\Delta$, then
$\alpha^{-1}\circ\theta:\Delta\to{\rm Sym}\,\Omega$ corresponds to
a hypermap ${\mathcal H}^{\alpha}$. Therefore, an \emph{operation}
on hypermaps is any transformation of hypermaps induced by a group
automorphism of $\Delta$. The hypervertices, hyperedges and
hyperfaces of ${\mathcal H}^{\alpha}$ are, respectively, the
orbits of $\langle r_1^{\alpha}, r_2^{\alpha}\rangle$, the orbits
of $\langle r_2^{\alpha}, r_0^{\alpha}\rangle$ and the orbits of
$\langle r_0^{\alpha}, r_1^{\alpha}\rangle$ on $\Omega$. If
$\alpha$ is an inner automorphism then ${\mathcal
H}^{\alpha}\cong{\mathcal H}$ for all $\mathcal H$, so we have an
induced action of the outer automorphism group ${\rm Out}\,\Delta
= {\rm Aut}\,\Delta/{\rm Inn}\,\Delta$ as a group $\Phi$ of
operations on isomorphism classes of hypermaps. This action is
faithful, as shown by L. James \cite{Lynne}. The same can be said
about oriented hypermaps, with ${\rm Out}\,\Delta^+$ acting as a
group $\Phi^+$ of operations. L. James \cite{Lynne} also proved
that ${\rm Out}\,\Delta^+ \cong GL(2,\mathbb{Z}_2)$, a very
important result for the classification of all operations on
oriented hypermaps (see \cite{GarethDaniel} for details).

\section{Algebraic Hypermaps}

$\Delta$ and $\Delta^+$ are, respectively, the full automorphism
group and the orientation-preserving automorphism group of the
universal hypermap ${\cal{\widetilde{H}}}$. This hypermap is
called \emph{universal} because any hypermap is the quotient of
${\cal{\widetilde{H}}}$ by some subgroup $H \leq \Delta$, known as
the \emph{hypermap subgroup} (which is unique up to conjugacy). If
$H \unlhd \Delta$, we say that the hypermap is \emph{regular}
since, when this occurs, the hypermap has the highest possible
number of symmetries. A regular hypermap can be represented,
algebraically, by a four-tuple ${\cal H}=(\Delta/H, h_0,h_1,h_2)$
where $h_0^2=h_1^2=h_2^2=1$ and $\langle h_0, h_1, h_2 \rangle =
\Delta /H$, the \emph{monodromy group}, $Mon ({\cal H})$, of the
hypermap. Similarly, an oriented regular hypermap can be regarded
as a triple ${\cal H}^+=(\Delta^+/H,x,y)$, with $\Delta^+/H=
\langle x, y \rangle$ being the monodromy group of the oriented
regular hypermap. From a topologically point of view, $x$ can be
interpreted as the permutation that cyclic permutes the hyperdarts
(oriented hyperedges) based on the same hypervertex, and $y$ the
permutation that cyclic permutes the hyperdarts based on the same
hyperface, according to the chosen orientation. An oriented
hypermap is called \emph{chiral} if it is not invariant under the
operation that reflects the oriented hypermap, inverting its
orientation (which is the same as saying that ${\cal H}^+$ admits
no orientation-reversing automorphism). If $(xy)^2=1$, ${\cal
H}^+$ is a map.

\section{The Duality Group} \label{thedualitygroup}

Our aim is to study what we will call the \emph{duality
group}\index{duality!group} of a hypermap. Some work has been done
on chirality groups \cite{ChiralityGroup} and there is no reason
not to extend that notion to duality or other hypermap operations.
These operations, as we have mentioned before, come from outer
automorphisms of $\Delta$ and by choosing the right group
$\Delta^*$, containing $\Delta$, we can look at duality as the
result of sending a hypermap subgroup to its conjugate in
$\Delta^*$. To build this group, we should add an element $t$, of
order 2, transposing $r_0$ and $r_2$ and fixing $r_1$. Hence, we
can define $\Delta^*$ in the following way:
$$\Delta^*=\Delta \rtimes C_2=\langle r_0,r_1,r_2,t: r_i^2=t^2=1,
r_0^t=r_2, r_1^t=r_1\rangle $$ This also means that $\Delta$ is a
normal subgroup of index 2 of $\Delta^*$. Therefore, each
conjugacy class of subgroups $H \leq \Delta$ is either a
$\Delta^*$-conjugacy class (if the hypermap $\cal H$ is self-dual,
which occurs when ${\cal H}\cong {\cal H}^t$) or paired with
another $\Delta$-conjugacy class, containing $H^t$ (if the
hypermap $\cal H$ is not self-dual). This last observation is a
general one and it is true for every kind of hypermap. However, we
will only deal with regular hypermaps and these have normal
subgroups as hypermap subgroups, which means that $H$ is conjugate
only to itself in $\Delta$. So, if a hypermap is self-dual, the
group $H$ is invariant under that specific outer automorphism of
$\Delta$ (conjugation in $\Delta^*$).

\begin{theorem}
Let $N$ be a normal subgroup of $\Delta$ and let $G=\Delta/N$.
Then the following are equivalent:

\begin{itemize}
\item[i)] $N^t=N$

\item[ii)] $N$ is normal in $\Delta^*$ \quad $ \hfill \square$
\end{itemize}

\end{theorem}

\noindent Because
$$\Delta^*=\langle r_0,r_1,r_2,t: r_i^2=t^2=1,
r_0^t=r_2, r_1^t=r_1\rangle = $$ $$=\langle r_0,r_1,t:
r_0^2=r_1^t=t^2=1, \mbox{ } r_1^t=r_1\rangle= \langle r_1,t
\rangle
* \langle r_0 \rangle \cong V_4*C_2 \cong \Gamma$$ we can build a
functor from hypermaps ($H \leq \Delta$) to maps ($H \leq \Delta^*
\cong \Gamma$) and, depending on the chosen isomorphism between
$\Delta^*$ and $\Gamma$, this is the Walsh functor \cite{Walsh},
representing a hypermap as a bipartite map or one of its duals.

If $\cal H$ is a regular hypermap with hypermap subgroup $H$ then
$H$ is normal in $\Delta$. The largest normal subgroup of
$\Delta^*$ contained in $H$ is the group $H_{\Delta}=H \cap H^t$
and the smallest normal subgroup of $\Delta^*$ containing $H$ is
the group $H^{\Delta}=HH^t$. These correspond, respectively, to
the smallest self-dual hypermap that covers $\cal H$, and the
largest self-dual hypermap that is covered by $\cal H$.

\begin{figure}[h!]
\begin{center}
\includegraphics [scale=0.7]{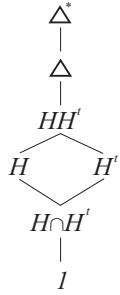}
\caption{$H_{\Delta}$ and $H^{\Delta}$.}
\end{center}
\end{figure} Like chirality, the duality operation is an operation of order
2 (see \cite{GarethDaniel} for the complete classification of
hypermap operations of finite order) and some of the results that
were stated for chirality and chirality groups
\cite{ChiralityGroup} also work here, with similar proves.
Whenever this is case, we will not give a demonstration of the
result since the reader can easily adapt the one available in
\cite{ChiralityGroup}. The following proposition is a good example
of what we have just mentioned:

\begin{proposition}\label{dualityproposition}
The groups $H^{\Delta}/H$, $H/H_{\Delta}$, $H^{\Delta}/H^t$ and
$H^t/H_{\Delta}$ are all isomorphic to each other. $\hfill
\square$
\end{proposition}

 \vskip 0.3cm

This common group will be called the \emph{duality group}
\index{duality!group}$D({\cal H})$ of $\cal H$ and its order the
\emph{duality index}\index{duality!index} $d$ of $\cal H$. The
duality index is somehow a way to measure how far the hypermap is
from being self-dual. If the duality index is 1 then the hypermap
is self-dual; and the bigger that index, the more distant the
hypermap is from being self-dual.

\begin{proposition}
The duality group $D({\cal H})$ of a regular hypermap $\cal H$ is
isomorphic to a normal subgroup of the monodromy group $Mon({\cal
H})$.
\end{proposition}
\noindent \emph{Proof:} The same as in \cite{ChiralityGroup}, with
duality group instead of chirality group.    \hfill$\square$

\vskip 0.3cm

Then, another possible way to understand the duality group is to
look at it as the minimal subgroup $D({\cal H}) \trianglelefteq
Mon({\cal H})$ such that ${\cal H}/D({\cal H})$ is a self-dual
hypermap. If $D({\cal H}) = Mon({\cal H})$ or, equivalently,
$H^{\Delta}=\Delta$ we say that the hypermap has \emph{extreme
duality index}.

\bigskip

We have extended $\Delta$ to $\Delta^*$ by adjoining $t$ such that
$t^2=1, \quad r_0^t=r_2, \quad r_2^t=r_0, \quad r_1^t=r_1.$ Then,
$$x^t=r_2r_1=y^{-1}, \quad y^t=r_1r_0=x^{-1}.$$
We will denote this kind of duality on oriented regular hypermaps
by $\beta-duality$ (chiral-duality)\index{chiral-duality}. On the
other hand, conjugation by $r_1t$ induces $$x \mapsto y, \quad y
\mapsto x,$$interchanging generators. This will be called
$\alpha-duality$ (orientation-preserving duality). The
relationship between the two ($\alpha \mbox{ and } \beta$) will be
dealt briefly at the end of this paper. From now on, to simplify
the writing, whenever we refer to \emph{duality} we mean
$\alpha-duality$, the one that preserves orientation. An hypermap
is called \emph{self-dual} if it is invariant under this duality
operation.

\bigskip

From an orientable hypermap we can choose two possible oriented
hypermaps. If ${\cal H}=(G,r_0, r_1, r_2)$, let ${\cal
H}^+=(G^+,x,y)$ be one of the oriented hypermaps associated with
$\cal H$. The duality group $D({\cal H}^+)$ is the minimal normal
subgroup of $G^+$ such that ${\cal H}^+ /D({\cal H}^+)$ is a
self-dual hypermap. It follows that $D({\cal H}^+) \unlhd G^+$ and
we say that ${\cal H}$ has extreme duality index if $D({\cal
H}^+)=G^+$. (Hence, a hypermap has extreme duality index
 if its duality group is equal to its
monodromy group).

\section{Duality index}

We can now easily prove the following theorem:

\begin{theorem}\label{sdt}

For every $k \in \mathbb{N}$, there is a self-dual oriented
regular hypermap with order $k$.
\end{theorem}
\noindent \emph{Proof}: Let $G$ be the cyclic group of order $k$
generated by $g$. If we take $G=\langle g\rangle $ and ${\cal
H}=(G,g,g)$, the hypermap with monodromy group $G$, then there is
an automorphism of $\cal H$ that interchanges the two generators
(they are both equal to $g$, in this case). Hence, the hypermap is
self-dual. \hfill$\square$

\bigskip

But this last result also means that for every $k \in \mathbb{N}$
there is a hypermap ${\cal H}=(G,a,b)$ such that $|G|/d=k$, with
$d$ being the duality index of $\cal H$. We just have to take
$G=\langle g \rangle$, as the cyclic group of order $k$, and the
hypermap ${\cal H}=(G,g,g)$, as in the previous proof. Because
$(G,g,g)$ is self-dual, $d=1$ and we have $|G|/d=|G|=k$. From now
on, we will call $|G|/d$ the \emph{duality coindex} of a hypermap
of monodromy group $G$.

\bigskip

Can we prove a similar theorem as Theorem \ref{sdt} using only
hypermaps that are not self-dual (for which $d \neq 1$)? The proof
we provide bellow will give the reader not only an affirmative
answer but also the presentation of the monodromy groups of those
hypermaps.

\begin{theorem}\label{quotient}
If $k \in \mathbb{N}$, there is a non self-dual  oriented regular
hypermap ${\cal H}=(G,a,b)$ with duality coindex
$k$\index{duality!coindex}.
\end{theorem}
\noindent \emph{Proof}:  Given $k \geq 3$, we can choose, by
Dirichlet's Theorem, a prime $q \equiv 1 \mbox{ mod }(k)$. Let
$G=\langle g,h|h^q=1, g^k=1, h^g=h^u \rangle \cong C_q \rtimes
C_k$, where $u \in \mathbb{Z}_q$ has multiplicative order k, $C_q
=\langle h \rangle$ and $C_k =\langle g \rangle$. Then, if $h=ab$
and $g=a$, we have:
$$ G= \langle a,b|(ab)^2=a^k=1, (ab)^a=(ab)^u \rangle $$
The duality group of this hypermap is the smallest normal subgroup
$N$ of $G$ such that the assignment $a\mapsto b$, $b \mapsto a$
induces an automorphism of $G/N$. We obtain this quotient by
adding extra relations, substituting $a$ for $b$ and $b$ for $a$
in the original ones.\footnote{This method will be used several
times in the next pages and the group $N$, in similar contexts,
will always mean the smallest normal subgroup $N$ of the monodromy
group $G$ such that the interchange of generators induces an
automorphism of $G/N$.} In this case, we just have to add these
relations: $b^k=1$ and $(ba)^b=(ba)^u$.Hence:
$$G/N=\langle (ab)^2=a^k=b^k=1, (ab)^a=(ab)^u, (ba)^b=(ba)^u \rangle$$
But $(ab)^a=ba$, so $ba=(ab)^u$, $ab=(ba)^u$. It follows that
$ab=(ab)^{u^2}$ or, equivalently, $(ab)^{u^2-1}=1.$ Because $k
\geq 3$, we have $u \neq \pm 1 \mbox{ mod } q \Rightarrow u^2-1
\neq 0\mbox{ mod } q \Rightarrow (u^2-1,q)=1.$ So,
$$(ab)^q=(ab)^{u^2-1}=1 \Rightarrow ab=1 \Rightarrow  b=a^{-1}.$$
Thus $G/N=\langle a|a^k=1 \rangle \cong C_k$. Therefore $|G/N|=k$
and, since $G$ is not cyclic, the hypermap ${\cal H}=(G,a,b)$ is
not self-dual.

\medskip

If $k=2$ we take $G=C_6$ generated by the pair $(1,4)$, with
presentation:
$$C_6=\langle x,y|x^6 =1, x^4=y\rangle .$$
\noindent Hence, considering $N$ as the duality group,
$$|C_6/N|=|\langle x,y|x^6=1, x^4=y, y^6=1, y^4=x\rangle |=3.$$
Then, $(C_6,x,y)$ has duality index $\frac{6}{3}=2$.

\medskip

For $k=1$, all we have to do is to choose any hypermap with
extreme duality index. \hfill$\square$
\bigskip

Now, another question can be asked: for each $d \in \mathbb{N}$,
is it possible to find at least one hypermap with that duality
index? And can we make some restrictions in the available
hypermaps we are allowed to choose? The first question is not
difficult to be answered:

\bigskip

\begin{theorem} \label{dualindextheorem}
For every $d \in \mathbb{N}$, there is an  oriented regular
hypermap with duality index equal to $d$.
\end{theorem}
\emph{Proof}: Let $G$ be the cyclic group of order $d$ generated
by $g$. If we take $G=\langle g\rangle $ and ${\cal H}=(G,g,1)$,
the hypermap with monodromy group $G$, then its duality group must
be equal to $G$, which means that the hypermap has an extreme
duality index $|G|=d$. $\hfill \square$

\bigskip
\noindent \textbf{Remarks}: a) Obviously, ${\cal H}=(G,1,g)$ also
works here. In fact, for any duality index, we can always find,
not just one, but two hypermaps with that extreme duality index
(which is not surprising since these two hypermaps are duals of
each other). b) It follows from the proof of this last theorem
that for every $n \geq 1$ there is an oriented regular hypermap
with cyclic duality group (the monodromy group of the hypermap
${\cal H}=(C_n,g,1)$ with $C_n=\langle g\rangle $).

\bigskip

It is now clear that we can get any duality index using hypermaps
that have extreme duality index. Can we achieve the same result
only with hypermaps that do not have extreme duality index? Before
we answer that question, we need to introduce some results and
definitions about direct products of hypermaps.

\section{Direct Products and Duality groups}

If $\cal H$ and $\cal K$ are oriented regular hypermaps with
hypermaps subgroups $H$ and $K \leq \Delta^+$, respectively, then:

\begin{definition}
The \emph{least common cover} ${\cal H}\vee {\cal K}$ and the
\emph{greatest common quotient} ${\cal H} \wedge {\cal K}$ are the
oriented regular hypermaps with hypermap subgroups $H\cap K$ and
$\langle H,K\rangle =HK$ respectively.
\end{definition}

If ${\cal H}=(D_1,R_1,L_1)$ and ${\cal K}=(D_2,R_2,L_2)$ let
$D=D_1 \times D_2$ and the permutations $R$ and $L$ be the ones
that act on $D$ induced by the actions $\rho \mapsto R_i$,
$\lambda \mapsto L_i$ of $\Delta^+$ on $D_1$ and $D_2$. If this
action is transitive on $D$, we call ${\cal H}\times {\cal
K}=(D,R,L)$, the \emph{oriented direct product} of $\cal H$ and
$\cal K$ with hypermap subgroup $H\cap K$.

\begin{lemma} \emph{\cite{BJ}}
If $\cal H$ and $\cal K$ are oriented regular hypermaps, then the
following conditions are equivalent:
\begin{itemize}
\item[i)] $\Delta^+$ acts transitively on $D$; \item[ii)] ${\cal
H}\wedge {\cal K}$ is the oriented hypermap, with one dart;
\item[iii)] $HK=\Delta^+$.\hfill$\square$
\end{itemize}
\end{lemma}

If these conditions are satisfied we say that $\cal H$ and $\cal
K$ are \emph{oriented orthogonal} and we use the notation ${\cal
H} \perp {\cal K}$. Then, ${\cal H} \times {\cal K}$ is well
defined and isomorphic to ${\cal H}\vee {\cal K}$ with monodromy
group $ Mon({\cal H} \times {\cal K})=Mon({\cal H})\times
Mon({\cal K})$. Having in mind that $\cal H$ has extreme duality
index if and only if $HH^d=\Delta^+$, we have, as an important
example, the following result:

\begin{lemma}
$\cal H$ has extreme duality index $\Leftrightarrow {\cal H} \perp
{\cal H}^d$. \hfill$\square$
\end{lemma}

Once again, we can adapt one of the theorems for chirality groups
 \cite{ChiralityGroup}, writing it in this new context of duality:

\begin{theorem} \label{producthypermap}
Let ${\cal H}$ and ${\cal K}$ be oriented regular hypermaps, with
hypermap subgroups $H$ and $K$, such that $\cal K$ has extreme
duality index and covers $\cal H$. Then the product ${\cal
L}={\cal K}\times {\cal H}^d$ is an oriented regular hypermap with
duality group $D({\cal L})\cong H/K$.
\end{theorem}

\noindent \emph{Proof:} The same as in \cite{ChiralityGroup},
substituting chirality by duality. \hfill$\square$

\medskip
We can now answer the question we have raised at the end of the
previous section:

\begin{theorem}\label{anyd}
For every $d \in \mathbb{N}$ there is an oriented regular hypermap
with non extreme duality index $d$.
\end{theorem}
\noindent \emph{Proof}: Let $K$ be a normal subgroup of $\Delta^+$
such that $\Delta^+/K=C_{2d}$. Then ${\cal K}=(C_{2d},g,1)$, with
$C_{2d}=\langle g\rangle $, is orientably regular and has extreme
duality index. If we take $H$ such that $K \leq H$ and $|H:K|=d$
then $|\Delta^+:H|=2$, which means that $H \trianglelefteq
\Delta^+$ and $\cal H$ is orientably regular. Hence, the hypermap
${\cal L}={\cal K}\times {\cal H}^t$ is orientably regular and
$$|Mon( {\cal L})|= |Mon ({\cal K})|\cdot|Mon({\cal H}^t)|=2d \times
2=4d.$$ Then, by Theorem \ref{producthypermap}:
$${\cal D}({\cal L})=H/K,$$
and $|H:K|=d$. $\cal L$ does not  have extreme duality index
because $|Mon ({\cal L})|=4d>d$. \hfill$\square$

\bigskip

This is not only true for \emph{hypermaps} but also for
\emph{maps}:

\begin{theorem}\label{dd}
For every $d \in \mathbb{N}$ there is an oriented regular map with
(non extreme) duality index equal to $d$.
\end{theorem}
\noindent \emph{Proof}: Let $D_{2m}=\langle
x,y|x^m=y^2=(xy)^2=1\rangle $ be the dihedral group of order $2m$.
If we take ${\cal M}=(D_{4d},x,y)$, then, considering $N$ as
before, we will have: $D_{4d}/N \cong D_4$. Therefore, $|N|=4d/4
=d$. $\hfill \square$

\medskip

Although a map is a special case of a hypermap (when $(xy)^2=1$),
Theorem \ref{dd} is not a stronger version of Theorem \ref{anyd},
since Theorem \ref{anyd} allows us to get not just hypermaps but
\emph{proper} hypermaps (hypermaps that are not maps), which is
also an important restriction.

\medskip

A group is called \emph{strongly self-dual} if for all its
generating pairs there is an automorphism of $G$ interchanging
them. A good example of one of these groups is the quaternion
group. In the next section, we will use a generalization of the
quaternion group to find infinite families of proper hypermaps
with non extreme duality indexes.

\section{Generalized quaternion groups}

\begin{definition}
If $w=e^{i \pi n} \in \mathbb{C}$, the matrices:
$$ x=\left(
\begin{array}{cc}
  w & 0 \\
  0 & \overline{w} \\
\end{array}
\right) , y=\left(
\begin{array}{cc}
  0 & 1 \\
  -1 & 0 \\
\end{array}
\right)
$$
generate a subgroup $Q_{2n}$ of order $4n$ in $GL(2,\mathbb{C})$
with presentation \cite{Johnson}: $$ \langle x,y|x^n=y^2,
x^{2n}=1, y^{-1}xy=x^{-1} \rangle$$ which is called the
\emph{generalized quaternion group}\index{generalized quaternion
group}.
\end{definition}

As we have proved in Theorem \ref{anyd}, we can have a hypermap of
any non extreme duality index. However, that proof does not show
us the presentation of the monodromy group of any of those
hypermaps. Various explicit examples can be obtained  using
generalized quaternion groups.

\begin{theorem}\label{odd}
If $d$ is odd or $d \equiv 0 \mbox{ }(mod \mbox{ }4)$, there is an
oriented regular hypermap with generalized quaternion monodromy
group, which has a non extreme duality index equal to d.
\end{theorem}
\noindent \emph{Proof}: a) $d$ is odd:

Let $n=2+4k$, $k=0,1,2,...$ . If we take $G$ to be the generalized
quaternion group of order $4n$ then $|G|=8+16k$ and has
presentation:
$$
G=\langle x,y|x^{2+4k}=y^2, x^{4+8k}=1, y^{-1}xy=x^{-1}\rangle.$$
If we take $N$ to be the smallest normal subgroup of $G$ such that
the assignment that interchanges the two generators induces an
automorphism then $G/N$ (which is obtain from $G$ adding new
relations) is the quaternion group and has order 8. But $
|N|=|G|/|G/N|$. Hence, $|N|=8+16k/8=2k+1 (\mbox{for } k=0,1,...)$
From this, we can conclude that for $d$ odd there is a hypermap
with monodromy group $G$ and a non extreme duality index (since
$|G/N|=8 \neq 1$) equal to $d=2k+1$.

\medskip
b) $d \equiv 0 \mbox{ }(mod \mbox{ }4)$:

Let $n=4k$, $k=1,2,...$ If we take $G$ to be the generalized
quaternion group of order $4n$ then $|G|=16k$ and has
presentation:
$$
G=\langle x,y|x^{4k}=y^2, x^{8k}=1, y^{-1}xy=x^{-1}\rangle $$ If
we take $N$ to be the smallest normal subgroup of $G$ such that
the assignment that interchanges the two generators induces an
automorphism then
$$G/N=\langle x,y|x^{4k}=y^2, x^{8k}=1,
y^{-1}xy=x^{-1}, y^{4k}=x^2, y^{8k}=1, x^{-1}yx=y^{-1} \rangle$$

Using the third and sixth relations, we have
$(y^{-1}xy)yx=x^{-1}(xy^{-1})=y^{-1}.$ Therefore, applying the
first relation $y^2=x^{4k}$, we have: $y^{-1}xx^{4k}x=y^{-1}
\Rightarrow x^{4k+2}=1.$ Then, using the second relation:
$x^{4k+2}=x^{8k}\Rightarrow x^{4k-2}=1$. Hence, $x^{4k+2}=x^{4k-2}
= 1 \Rightarrow x^4=1.$

From the first relation $x^{4k}=y^2$ we can now conclude that
$y^2=1$ and, from the fourth one, that $x^2=1$. Therefore, the
presentation of the group G reduces to $\langle x,y |
x^2=y^2=(xy)^2=1 \rangle$, which defines a Klein 4-group.
 We have proved, this way, that $G/N$ has order 4. Hence $$|N|=16k/4=4k \quad ,\quad k=1,2,...$$

This means that, for $d \equiv 0 \mbox{ }(mod \mbox{ }4)$, there
is a hypermap with monodromy group $G$ has a non extreme duality
index (since $|G/N|=4 \neq 1$) equal to $d=4k$. $\hfill \square$
\begin{corollary}
Every cyclic group of odd order or of order multiple of 4 can be a
duality group of an oriented regular hypermap with non extreme
duality index and generalized quaternion monodromy group.
\end{corollary}
\noindent \emph{Proof}: In the previous proof $N=\langle x^4
\rangle \cong C_{1+2k}$, in a); and $N=\langle x^2 \rangle \cong
C_{4k}$, in b). $\hfill \square$

\bigskip

In the proof of the Theorem \ref{odd}, $G/N$ is the quaternion
group and any hypermap which has that group as monodromy group is
self-dual. But all generating pairs are equivalent under
automorphisms of the quaternion group. Then, there is only one
(self-dual) hypermap, up to isomorphism, with monodromy group
being the quaternion group.

\begin{theorem}
Let $n$ be odd. Then, the generalized quaternion group $$G=
\langle x,y|x^n=y^2, x^{2n}=1, y^{-1}xy=x^{-1}\rangle $$ of order
$4n$ is the monodromy group of an oriented regular hypermap with
extreme duality index.
\end{theorem}
\noindent \emph{Proof}: If we take $N$ to be the smallest normal
subgroup of $G$ such that the assignment that interchanges the two
generators induces an automorphism then $$G/N=\langle x,y|x^n=y^2,
x^{2n}=1, y^{-1}xy=x^{-1}, y^n=x^2, y^{2n}=1,
x^{-1}yx=y^{-1}\rangle .$$ Hence, we have $x^{-1}yx=y^{-1}$ (last
relation) but also $x^{-1}=y^{-1}xy$ (third relation). Therefore,
$ y^{-1}xyyx=y^{-1} \Rightarrow y^{-1}xy^2x=y^{-1}.$ Using the
first relation in this last equality, we have $y^{-1}xx^nx=y^{-1}
\Rightarrow x^{n+2}=1.$ Let $k$ be the order of $x$. Then, since
$k|(n+2)$ and $n$ is odd, $k$ must also be odd. But from the
second relation we also know that $x^{2n}=1$ and, consequently,
$k|2n$. Therefore, $k|n$. If odd $k$ divides $n$ and $n+2$, then
$k=1$ (and we have $x=1$). Because $y^n=y^2$ and $y^2=x^n$, we
have $y^n=y^2=1$. Since $n$ is odd, $y=1$. Hence, $|G/N|={1}$,
which means that the hypermap has extreme duality index. $\hfill
\square$

\begin{corollary}
There are infinitely many  oriented regular hypermaps with extreme
duality index and generalized quaternion group as monodromy group.
\hfill$\square$
\end{corollary}

Every hypermap having the generalized quaternion group (with the
presentation given in our definition) as monodromy group has
chirality index equal to 1. This can easily be checked because if
we want to obtain a reflexible hypermap as a quotient of the
original one, we just have to add the following relations to the
ones that we already have for the generalized quaternion group:
$x^{-n}=y^{-2}$, $x^{-2n}=1$ and $yx^{-1}y^{-1}=x$ (substituting
$x$ by $x^{-1}$ and $y$ by $y^{-1}$ in the original relations) .
However, these relations do not change the presentation of the
group. Hence, all the theorems above (where the generalized
quaternion group appears in the proof) are, in fact, about
reflexible (non chiral) hypermaps.

\section{Chiral duality}
As we have previously noticed, there are two types of duality
induced by the following automorphisms of $\Delta^+$:

$$\alpha : x \mapsto y; \quad y \mapsto x,$$
$$\beta : x \mapsto y^{-1}; \quad y
\mapsto x^{-1}.$$ Since the automorphisms of $\Delta^+$ which
induce them are conjugate in $Aut(\Delta^+)$, both dualities have
the same general properties (the groups which arise as
$\alpha-duality$ groups are the same that arise as $\beta$-duality
groups \cite{GarethDaniel}). Nevertheless, their effect on a
\emph{specific} hypermap might be distinct. To make this
observation clear to the reader, we will give some examples of
hypermaps such that:
\begin{itemize}
\item[a)] $|{\cal D}_{\alpha}({\cal H})| \neq |{\cal
D}_{\beta}({\cal H})|$

\item[b)] ${\cal D}_{\alpha}({\cal H}) \cong {\cal
D}_{\beta}({\cal H})$
\end{itemize}

\textbf{Examples}

\bigskip

a) We can take ${\cal H}=(G,x,y)$ with $G=\langle x,y|x^4=y^4=1,
\quad xy=y^2x^2 \rangle.$ $|G|=20 \mbox{ (this order can easily be
checked using GAP \cite{GAP})}$. Then
$$ G/N_{\alpha} =\langle x,y|x^4=y^4=1,
\quad xy=y^2x^2, \quad yx=x^2y^2 \rangle.$$ Using the two last
relations, we have: $xyyx=y^2x^2x^2y^2 \Leftrightarrow xy^2x=1
\Leftrightarrow y^2=x^2.$ Therefore, $G/N=\langle x | x^4=1
\rangle$ and $|G/N_{\alpha}|=4.$ However: $$ G/N_{\beta} =\langle
x,y|x^4=y^4=1, \quad xy=y^2x^2, \quad y^{-1}x^{-1}=x^{-2}y^{-2}
\rangle=G.$$ Hence $|G/N_{\beta}|=20$. It follows that $G$ is
$\beta$-self-dual but not $\alpha$-self-dual.

\bigskip

b) If ${\cal H}=(G,x,y)=(A_5,(12345),(123))$ then
$D_{\alpha}({\cal H}) \cong A_5$ because the hypermap has extreme
$\alpha$-duality index. But ${\cal
H}^{\beta}=(G,y^{-1},x^{-1})=(A_5,(132),(15432))$. Hence, we still
have two permutations of different order. This means that the
hypermap cannot be $\beta$-self-dual and, because $A_5$ is simple,
we can conclude that it must have extreme $\beta$-duality index.
It follows that $D_{\beta}({\cal H}) \cong D_{\alpha}({\cal H})
 \cong A_5$.

\bigskip

\noindent \emph{Acknowledgments}: I would like to thank Professor
Gareth Jones for the useful discussions about this work.

\end{document}